\documentclass[11pt,leqno]{article}

\usepackage[english]{babel}

\usepackage{amsmath}

\usepackage{amssymb}

\usepackage{amsfonts}

\usepackage{mathabx}

\numberwithin{equation}{section}

\newcommand{\C}{\mathcal{C}}

\newcommand{\N}{\mathcal{N}}

\renewcommand{\S}{\mathcal{S}}

\newcommand{\T}{\mathcal{T}}

\newcommand{\U}{\mathcal{U}}

\renewcommand{\mod}{\mathrm{Mod}}

\newcommand{\R}{\mathbb{R}}

\newcommand{\OO}{\mathcal{O}}

\newcommand{\G}{\mathcal{G}}

\newcommand{\tho}{\mathit{T}\mathcal{H}\mathit{om}}

\newcommand{\wtens}{\overset{\mathrm{w}}{\otimes}}

\renewcommand{\dim}{\textbf{Proof.}}

\newcommand{\qed}{\nopagebreak \phantom{} \hfill $\Box$ \\}

\newtheorem{teo}{Theorem}[section]

\newtheorem{df}[teo]{Definition}

\newtheorem{cor}[teo]{Corollary}

\newtheorem{oss}[teo]{Remark}

\newtheorem{lem}[teo]{Lemma}

\author{Luca Prelli}

\title{De Rham theorem for Schwartz functions on Nash manifolds.}

\date{}

\begin{document}

\maketitle

\begin{abstract}
In \cite{AG10} the authors proved the de Rham theorem for Schwartz functions on affine Nash manifolds. Here we simplify the proof and generalize their result to the case of non affine Nash manifolds.
\end{abstract}

\section{Introduction}

In \cite{AG10} the authors 
proved analogs of the de Rham theorem on affine Nash manifolds for de Rham complexes
with coefficients in Schwartz functions and generalized Schwartz
functions. Using that they proved a version of the Shapiro lemma for
Schwartz functions on affine Nash manifolds.\\

In this paper we generalize their de Rham theorem to the case of non affine Nash manifolds, using different techniques.
Let $M$ be a Nash manifold. The de Rham theorem states that the cohomology of the de Rham complex with Schwartz coefficients $DR_\S(M)$ is isomorphic to the compact support cohomology of $M$.
Instead of considering the de Rham complex with Schwartz coefficients we consider the de Rham complex with coefficients
in generalized Schwartz functions $DR_\G(M)$ and prove that its cohomology is isomorphic to the cohomology of $M$. Then we obtain the result by duality.
In order to do that we need some sheaf theory and the notion of restricted topology on $M$, roughly speaking, we only consider semi-algebraic open subsets and finite covers. With this topology, generalized Schwartz functions are acyclic with respect to global sections and $DR_\G(M)$ is quasi-isomorphic to the complex $R\Gamma(M;\R_M)$ of global sections of the constant sheaf. 
This is known to be quasi-isomorphic to $DR_{-\infty}(M)$, the de Rham complex with coefficients in generalized functions.
By duality we obtain that $DR_\S(M)$ is quasi-isomorphic to $DR_c(M)$, the de Rham complex with compact support and its cohomology is isomorphic to the compact support cohomology of $M$.\\

\noindent \textbf{Acknowledgments.} We would like to thank D. Gourevitch, N. Honda and P. Polesello for their comments which helped
us to improve the presentation of this work.

\section{Restricted topologies and sheaves.}

Let $M$ be a topological space. We denote by $\mod(\R)$ the category of $\R$-vector spaces and by $\mod(\R_M)$ that of sheaves of $\R$-vector spaces on $M$. For classical sheaf theory we refer to \cite{KS90}. The following results are extracted from \cite{EP} and \cite{KS01}. Let us consider a family $\T$ of open subsets of $M$ satisfying:
\begin{equation*}\label{hytau}{ (\hspace{-0.5mm}\bigvarstar\hspace{-0.5mm})}  \
  \begin{cases}
    \text{(i) $\T$ is closed under finite unions and intersections and $\varnothing \in \T$},\\
    \text{(ii) any $U \in \T$ has finitely many connected components $\in \T$},\\
    \text{(iii) $\T$ is a basis for the topology of $M$}. \\
  \end{cases}
\end{equation*}


We may associate to $\T$ a Gro\-then\-dieck topology, called the restricted topology, in the following way:
a family $\U=\{U_i\}_i$ in $\T$ is a covering of $U \in \T$ if it admits a finite subcover. We denote by $\widetilde{M}$ the associated site. A sheaf on $\widetilde{M}$ is a functor (presheaf) $F:\T^{op} \to \mod(\R)$ satisfying the following gluing conditions: for each $U \in \T$ and for each covering $\U$ of $U$ in $\T$ the sequence
$$
0 \to F(U) \to \prod_{\U_i \in \U}F(U_i) \rightrightarrows \prod_{U_j,U_k \in \U}F(U_j \cap U_k)
$$
is exact. We call $\mod(\R_{\widetilde{M}})$ the category of sheaves of $\R$-vector spaces on $\widetilde{M}$. \\

Let  $\rho: M \to \widetilde{M}$ be the natural morphism of sites. The direct image $\rho_*:\mod(\R_M) \to \mod(\R_{\widetilde{M}})$ is defined by $\rho_*F(U)=F(U)$. It is a fully faithful functor.


\begin{df} A sheaf $F$ on $\widetilde{M}$ is $\T$-flabby if the restriction $F(U) \to F(V)$ is surjective for each $U,V \in \T$ with $U \supseteq V$. \\
\end{df}

The family of $\T$-flabby sheaves is injective with respect to the functor $\Gamma(\widetilde{M};\cdot)$.
Hence, for each $F \in \mod(\R_{\widetilde{M}})$, one may compute the cohomology of $R\Gamma(\widetilde{M};F)$, the right derived functor of $\Gamma(\widetilde{M};\cdot)$, by considering a resolution
$0 \to F_0 \to F_1 \to \cdots$ of $F$ and by taking the cohomology of the complex $\Gamma(\widetilde{M};F_\bullet)$.

\section{Nash manifolds.}

In this section we recall the definition of Nash manifolds, following \cite{BCR}, \cite{Shi} and \cite{AG08}.\\

Recall that a semi-algebraic subset of $\R^n$ is a subset defined by a finite number of polynomial equations and inequalities, or any finite union of such sets. A function is said to be semi-algebraic if its graph is semi-algebraic.\\

 We shall need the following two properties of semi-algebraic subsets, for which we refer to \cite{BCR}: (i) an open semi-algebraic subset has a finite number of connected components, (ii) thanks to the triangulation theorem, an open semi-algebraic subset has a finite cover consisting of semi-algebraic open subsets.

\begin{df} (i) A Nash function from an open semi-algebraic subset $U \subseteq \R^n$ to
an open semi-algebraic subset $V \subseteq \R^m$ is a semi-algebraic $\C^\infty$-function. The ring of $\R$-valued Nash functions on $U$ is denoted
by $\N(U)$. A Nash diffeomorphism is a Nash bijection whose inverse is
also Nash.\\

(ii) A semi-algebraic subset $M \subseteq \R^n$ is said to be a Nash submanifold
of $\R^n$ of dimension $d$ if, for every point $x \in M$, there exists a Nash
diffeomorphism $\phi:U \to U'$, where $U$ is an open semi-algebraic neighborhood of the origin in
$\R^n$ and $U'$ is an open semi-algebraic neighborhood of $x$ in $\R^n$, such that $\phi(0) = x$
and $\phi((\R^d \times \{0\}) \cap U) = M \cap U'$.
A Nash function from a Nash submanifold $M$ of $\R^m$ to a Nash
submanifold $N$ of $\R^n$ is a semi-algebraic $\C^\infty$-function.
\end{df}

In order to define Nash manifolds we need the following.

\begin{df} Let $M$ be a topological space and let $\T$ be a family of open subsets of $M$ satisfying 
$ { (\hspace{-0.5mm}\bigvarstar\hspace{-0.5mm})} $. A $\R$-space is a pair $(M,\OO_M)$
where $\OO_M$ is a subsheaf of $\R$-algebras over $\widetilde{M}$ of the
sheaf of real-valued functions.
A morphism between $\R$-spaces $(M,\OO_M)$ and $(N,\OO_N)$ is a morphism of sites $f:\widetilde{M}\to \widetilde{N}$, such that the induced morphism of sheaves
maps $\OO_N$ to $\OO_M$.
\end{df}

Take for $M$ a Nash submanifold of $\R^n$ and for $\T$ the family of all open subsets of $M$ which are semi-algebraic in $\R^n$. Then $\T$ satisfies ${ (\hspace{-0.5mm}\bigvarstar\hspace{-0.5mm})}$. 
For any  subset $U \in \T$ of $M$ we take the algebra $\N(U)$ of Nash functions. This defines a sheaf on $\widetilde{M}$.

\begin{df} (i) An affine Nash manifold is an $\R$-space which is isomorphic
as an $\R$-space to a closed Nash submanifold of $\R^n$. A morphism between two
affine Nash manifolds is a morphism of $\R$-spaces between them.\\

(ii) A Nash manifold is an $\R$-space $(M,\N_M)$ which has a finite
cover $(M_i)$ by open sets $M_i$ such that the $\R$-spaces $(M_i,\N_M|_{M_i}
)$ are isomorphic
as $\R$-spaces to affine Nash manifolds.
A morphism between Nash manifolds is a morphism of $\R$-spaces between
them.
\end{df}

By the open semi-algebraic subset of a Nash manifold we
mean its open subset in the  restricted topology. We shall need the following result (Remark I.5.12 of \cite{Shi}).

\begin{teo} \label{nash} Any Nash manifold can
be covered by a finite number of open submanifolds diffeomorphic (as Nash manifolds) to $\R^n$.
\end{teo}

We end the section with this lemma about the cohomology of the constant sheaf. Denote by $\R_M$ (resp. $\R_{\widetilde{M}}$) the constant sheaf on $M$ (resp. $\widetilde{M}$).

\begin{lem} \label{csh} 
(i) The adjunction morphism $\R_{\widetilde{M}} \to R\rho_*\R_M$ is an isomorphism.
(ii) In particular $R\Gamma(M;\R_M) \simeq R\Gamma(\widetilde{M};\R_{\widetilde{M}})$.
\end{lem}
\dim\ \ (i) Since the problem is local in $\widetilde{M}$, it is enough to prove the isomorphism on a finite cover of $M$, so by Theorem \ref{nash} we may assume $M=\R^n$. Recall that $\T$ is the family of open semi-algebraic subsets of $M$.
\begin{itemize}
\item[$\mathrm{(i)}_{\rm a}$] We first prove that $R^k\rho_*\R_M=0$ if $k \neq 0$. First note that, since $\rho_*$ sends injective objects to injective objects, we have $R\Gamma(U;R\rho_*\R_M) \simeq R\Gamma(U;\R_M)$ for each $U \in \T$. It is enough to prove that each $U \in \T$ has a finite cover $\{U_i\}$ in $\T$ such that $H^k(U_i;R\rho_*\R_M) \simeq H^k(U_i;\R_M)=0$ if $k \neq 0$. This follows from the fact that, thanks to the triangulation theorem, each $U \in\T$ has a finite open cover consisting of contractible semi-algebraic subsets and sections of the constant sheaf on contractible open sets have cohomology concentrated in degree zero (Corollary 2.7.7 of \cite{KS90}).
\item[$\mathrm{(i)}_{\rm b}$] It remains to prove that $\rho_*\R_M \simeq \R_{\widetilde{M}}$. 
    This follows from the fact that each open semi-algebraic subset has a finite number of connected components and when $U \in \T$ is connected $\Gamma(U;\rho_*\R_M) \simeq \R \simeq \Gamma(U;\R_{\widetilde{M}})$.
\end{itemize}
(ii) Since $\rho_*$ sends injective objects to injective objects we have $R\Gamma(M;\R_M) \simeq R\Gamma(\widetilde{M};R\rho_*\R_M) \simeq R\Gamma(\widetilde{M};\R_{\widetilde{M}})$. \qed \\

\section{De Rham complex with Schwartz coefficients.}

For the notions on topological vector spaces we shall need now, we refer to \cite{Gr}. A short exposition can be also found in the Appendix A of
\cite{CHM}.

\begin{df} We call a complex of topological vector spaces admissible if
all its differentials have closed images.
\end{df}
We will need the following classical facts from the theory of nuclear Fr\'echet (FN) and dual nuclear Fr\'echet (DFN)
spaces.

\begin{lem} \label{adm} Let
$
C_\bullet = 0 \to C_0 \to C_1 \to \cdots \to C_n \to 0
$
be a complex of FN spaces (or DFN spaces). Denote by $C_\bullet'$ the dual complex.
\begin{itemize}
\item[(i)] Suppose that $C_\bullet$ is admissible. Then the dual complex $C_{\bullet}'$ is also admissible and $H^k(C_\bullet)' \simeq H^{-k}(C_\bullet')$.
\item[(ii)] Suppose that $C_\bullet$ has
finite dimensional cohomology. Then $C_\bullet$ is admissible.
\end{itemize}
\end{lem}

In \cite{AG08} the authors defined, for any Nash manifold $M$, the FN space $\S(M)$
of Schwartz functions. When $M$ is a closed affine Nash submanifold of $\R^n$, 
a $\C^\infty$-function $s$ belongs to $\S(M)$
if  $Ps$ is bounded for any algebraic differential operator $P$ on
$M$.\\

Let $D_M$ denote the space of densities (with Nash coefficients) on $M$.
The space
of generalized Schwartz functions is defined by
$\G
(M) = (\S
(M) \otimes_{\N(M)} D_M)'$. In \cite{AG08} the authors proved that the correspondence $U \mapsto \G(U)$ defines a $\T$-flabby sheaf on $\widetilde{M}$. \\

Let $M$ be a Nash manifold. We denote by
\begin{itemize}
\item $DR_c(M)$: the de Rham complex of $M$ with coefficients in $\C^\infty_c(M)$, the space of $\C^\infty$-functions with compact support;
\item $DR_\S(M)$: the de Rham complex of $M$ with coefficients in $\S(M)$, the space of Schwartz functions;
\item $DR_{-\infty}(M)$: the de Rham complex of $M$ with coefficients in classical generalized functions, i.e. functionals
on compactly supported densities;
\item $DR_\G(M)$: the de Rham complex of $M$ with coefficients in $\G(M)$, the space of generalized Schwartz functions.
\end{itemize}

\begin{teo} \label{SDR} (i) The complexes $DR_{-\infty}(M)$ and $DR_\G(M)$ are quasi-isomorphic.

(ii) The complexes $DR_c(M)$ and $DR_\S(M)$ are quasi-isomorphic.
\end{teo}
\dim\ \ (i) 
The de Rham complex with coefficients in generalized functions is finite dimensional and quasi-isomorphic to $R\Gamma(M;\R_M)$. On the other side, since $\G$ is $\T$-flabby the de Rham complex with coefficients in generalized Schwartz functions is a resolution of $R\Gamma(\widetilde{M};\R_{\widetilde{M}})$.
We are reduced to prove $R\Gamma(M;\R_M) \simeq R\Gamma(\widetilde{M};\R_{\widetilde{M}})$ and the result follows from Lemma \ref{csh}. \\
\noindent (ii) Let us prove the assertion in several steps. For $\star=-\infty,\G$, we denote by $DR_\star(M)^\vee=DR_\star(M) \otimes_{\N(M)} D_M$.
\begin{itemize}
\item[$\mathrm{(ii)}_{\rm a}$] It is well known that the cohomologies of $DR_c(M)$ and $DR_{-\infty}(M)^\vee$ are finite dimensional and dual to each other.
\item[$\mathrm{(ii)}_{\rm b}$] The complex $DR_\S(M)$ is dual of $DR_\G(M)^\vee$ as a complex of DFN spaces, hence by Lemma \ref{adm} the cohomology groups of $DR_\S(M)$ are finite dimensional.
\item[$\mathrm{(ii)}_{\rm c}$] By (i) and $\mathrm{(ii)}_{\rm b}$, $DR_\S(M)$ is quasi-isomorphic to the dual of $DR_{-\infty}(M)^\vee$. Hence by $\mathrm{(ii)}_{\rm a}$ $DR_c(M)$ and $DR_\S(M)$ are quasi-isomorphic. \qed \\
\end{itemize}





As usual, we denote by $H^\bullet_c(M)$ (resp. $H^\bullet(M)$) the cohomology with compact support (resp. cohomology) of $M$.

\begin{cor} Let $M$ be a Nash manifold. Then
$
H^\bullet(DR_\S(M)) \simeq H_c^\bullet(M)$ and $H^\bullet(DR_\G(M)) \simeq H^\bullet(M).
$
In particular, they are finite dimensional.
\end{cor}
\dim\ \ It follows immediately from Theorem \ref{SDR} and the classical isomorphisms $H^\bullet(DR_c(M)) \simeq H_c^\bullet(M)$ and $H^\bullet(DR_{-\infty}(M)) \simeq H^\bullet(M)$.
\qed

\begin{oss} (i) After one proves the de Rham theorem for Nash manifolds, the relative
de Rham theorem and the Shapiro lemma will follow in the same way as
in \cite{AG10}, as pointed out by the authors.

(ii) One could also develop an analogue of the functors $\wtens$ and $\tho$ of \cite{KS96} in this setting and obtain a further generalization of the
de Rham theorem.
\end{oss}

\end{document}